\documentclass[reqno,11pt]{amsart}
\usepackage{amsthm, amsmath, amssymb,booktabs,xcolor,graphicx}
\usepackage[margin=1in]{geometry}

\usepackage[colorlinks]{hyperref}

\usepackage[noabbrev,capitalize]{cleveref}
\crefname{equation}{}{}

\allowdisplaybreaks

\usepackage{cancel}

\newtheorem{theorem}{Theorem}[section]
\newtheorem{proposition}[theorem]{Proposition}

\newtheorem{corollary}[theorem]{Corollary}

\newtheorem{conjecture}[theorem]{Conjecture}

\newtheorem*{question*}{Question}

\theoremstyle{definition}
\newtheorem{definition}[theorem]{Definition}

\newtheorem{question}[theorem]{Question}
\newtheorem*{definition*}{Definition}

\theoremstyle{remark}

\newcommand{\floor}[1]{\left\lfloor #1 \right \rfloor }

\newcommand{\paren}[1]{\left( #1 \right)}

\newcommand{\EE}{\mathbb{E}}
\newcommand{\CC}{\mathbb{C}}

\newcommand{\RR}{\mathbb{R}}

\newcommand{\NN}{\mathbb{N}}
\newcommand{\FF}{\mathbb{F}}
\newcommand{\ZZ}{\mathbb{Z}}

\title[Monochromatic solutions over the integers]{On monochromatic solutions to linear equations\\ over the integers}
\author[Dong]{Dingding Dong}
\author[Mani]{Nitya Mani}
\author[Pham]{Huy Tuan Pham}
\author[Tidor]{Jonathan Tidor}
\thanks{Mani was supported by the NSF Graduate Research Fellowship Program and the Hertz Graduate Fellowship. Pham was supported by a Clay Research Fellowship and a Stanford Science Fellowship. Tidor was supported by a Stanford Science Fellowship.}

\address{Department of Mathematics, Harvard University, Cambridge, MA 02138, USA}
\email{ddong@math.harvard.edu}

\address{Department of Mathematics, Massachusetts Institute of Technology, Cambridge, MA 02139, USA}
\email{nmani@mit.edu}

\address{School of Mathematics, Institute for Advanced Study, Princeton, NJ 08540, USA}
\email{htpham@caltech.edu}

\address{Department of Mathematics, Stanford University, Stanford, CA 94305, USA}
\email{jtidor@stanford.edu}

\begin{document}

\begin{abstract}
We study the number of monochromatic solutions to linear equations in a $2$-coloring of $\{1,\ldots,n\}$. We show that any nontrivial linear equation has a constant fraction of solutions that are monochromatic in any $2$-coloring of $\{1,\ldots,n\}$. We further study commonness of four-term equations and disprove a conjecture of Costello and Elvin by showing that, unlike over $\FF_p$, the four-term equation $x_1 + 2x_2 - x_3 - 2x_4 = 0$ is uncommon over $\{1,\ldots,n\}$.
\end{abstract}

\maketitle

\section{Introduction}\label{sec:intro}
In the 1990s, Graham, R\"odl and Ruci\'nski~\cite{GRR96} asked the following basic arithmetic question: suppose we 2-color the elements of $[n]:=\{1,\dots,n\}$ and count the number of monochromatic Schur triples, i.e., the number of triples $(x,y,z)\in[n]^3$ such that $x+y = z$ and $x,y,z$ receive the same color. How many such triples must be there? This question was solved by Robertson and Zeilberger~\cite{RZ98} (see also \cite{Dat03,Sch99}), who showed that the minimum number of monochromatic Schur triples in a 2-coloring of $[n]$ is asymptotically $n^2/11+O(n)$. 

More generally, given a $k$-term linear equation $a_1x_1+\dots+a_kx_k=0$ on variables $x_1,\dots,x_k$ with integer coefficients $a_1,\dots,a_k\in\ZZ\setminus\{0\}$, we can ask for the minimum number of monochromatic solutions to the equation in a 2-coloring of $[n]$. This problem has been studied for several linear equations (and systems of linear equations) including generalized Schur triples, $k$-term arithmetic progressions, and constellations \cite{BCG10,GKW24,LP12,PRS08,TW17}. While the answers to such questions are fairly well-understood when we consider 2-colorings over $\FF_q^n$, determining the minimum number of monochromatic solutions in 2-colorings of $[n]$ is far less understood.

In fact, the following basic question is still open: for which linear equations does a {uniformly random $2$-coloring} asymptotically minimize the number of monochromatic solutions? Informally, we say that a $k$-term linear equation is \textit{common over the integers} if any $2$-coloring of $[n]$ has at least as many monochromatic solutions asymptotically (as $n\to\infty$) as a uniformly random coloring. (We will give a formal definition in \cref{sec:prelim}.)
In 2006, Parrilo, Robertson, and Saracino \cite{PRS08} disproved a folklore conjecture by showing that three-term arithmetic progressions are {not common}, i.e.,  a random coloring does {not} minimize the number of monochromatic solutions to $x - 2y + z = 0$. The expected number of monochromatic solutions to this equation in a random 2-coloring of $[n]$ is $(1/8+o(1))n^2$. However, Parrilo, Robertson, and Saracino showed that the true minimum (call this number $T(n)$) satisfies
\[\paren{\frac{1675}{16384}+o(1)}n^2\leq T(n)\leq \paren{\frac{117}{1096}+o(1)}n^2.\]
In particular, $T(n)< (1/8+o(1))n^2$. The true value of $T(n)$ is still not known, even asymptotically.

More recently, Costello and Elvin~\cite{CE23} extended this result by proving that {all} $3$-term equations are uncommon over the integers. In their paper, they posed two conjectures related to this question of commonness of linear equations, both of which we resolve in this work. 

Over $\FF_p^n$, common equations are quite well-understood. Fox, Pham, Zhao \cite{FPZ21} completely characterized the common linear equations in finite field vector spaces. They proved that for any equation with non-zero coefficients, if the equation has odd length then it is common, and if the equation has even length then it is common if and only if its coefficients can be partitioned into pairs summing up to zero.

It is not hard to show that every uncommon equation over $\FF_p^n$ is uncommon over the integers as well (see also \cite[Lemma 2.1]{CE23}). Fox, Pham, and Zhao showed that if an equation is uncommon over $\FF_p^n$, then there exists a 2-coloring of $\FF_p$ with too few monochromatic solutions. Then for $m$ large, set $n=mp$ and 2-color $[n]$ by $m$ repeated copies of the given 2-coloring of $\FF_p$. For $m$ sufficiently large, this coloring will also have too few monochromatic solutions. Thus the work of Fox--Pham--Zhao implies that every even length equation whose coefficients cannot be partitioned into pairs summing up to zero is uncommon over the integers.

Note however that the work of Costello and Elvin shows that the integer setting and the finite field setting can have significantly different behaviors: all 3-term equations are common over $\FF_p^n$ but uncommon over the integers. Despite this, Costello and Elvin conjectured that the linear equation $a_1x_1+\dots+a_kx_k=0$ is common over the integers if and only if $k$ is even and $a_1,\dots,a_k$ can be partitioned into pairs summing up to zero \cite[Conjecture 4.1]{CE23}. As evidence for their conjecture, they proved that $x_1+x_2-x_3-x_4=0$ is common. We show that this conjecture does not hold by exhibiting the following counterexample.

\begin{theorem}\label{thm:even-uncommon}
    The linear equation $x_1+2x_2-x_3-2x_4=0$ is uncommon over the integers.
\end{theorem}

In the course of proving~\cref{thm:even-uncommon}, we will show that $ax_1+b x_2-ax_3-bx_4=0$ is common over the integers if and only if it is locally common over the integers (defined precisely in \cref{sec:prelim}). We further show that these notions are equivalent to a certain sequence of $n\times n$ matrices eventually being positive semidefinite. While falling short of a characterization of common $4$-term linear equations, the proof of~\cref{thm:even-uncommon} provides a recipe for checking if $ax_1+b x_2-ax_3-bx_4=0$ is common over the integers for any particular choice of $a$ and $b$.

In~\cite{CE23}, in addition to the upper bound implied by uncommonness, the authors additionally showed that for all $a,b\in\NN$, every 2-coloring of $[n]$ has $\Omega(n^2)$ monochromatic solutions to the equation $ax_1+ax_2-bx_3=0$. 

The authors conjectured that a similar result should hold for all linear equations; namely, for any linear equation $a_1x_1+\cdots+a_kx_k=0$ with at least one positive and one negative coefficient, every 2-coloring of $[n]$ should have at least $\Omega(n^{k-1})$ monochromatic solutions to the equation \cite[Conjecture 3.1]{CE23}. (The reason for the hypothesis on the signs of the coefficients is to guarantee that the equation has at least one solution in $\mathbb{N}$.) We verify that this conjecture is true.

\begin{theorem}\label{thm:lower}
For $k\geq3$ and $a_1, \ldots, a_k \in \ZZ\setminus\{0\}$, at least one of which is positive and one of which is negative, there exists $\delta >0$ such that for $n$ sufficiently large, every 2-coloring of $[n]$ has at least $\delta n^{k-1}$ monochromatic solutions to the equation $a_1x_1+\dots+a_kx_k=0$.
\end{theorem}

\cref{thm:lower} can also be generalized to larger systems of linear equations. 
\begin{theorem}\label{thm:lower-sys}
Let $B\in \mathbb{Z}^{w \times k}$ be an integer matrix of rank $w$. Let $r\geq 2$ be such that for $n$ sufficiently large, in every $r$-coloring of $[n]$ there is a monochromatic solution to the system of $w$ equations $Bx=0$. Then there exists $\delta=\delta(r,B)>0$ such that for $n$ sufficiently large, every $r$-coloring of $[n]$ has at least $\delta n^{k-w}$ monochromatic solutions to the system $Bx=0$.
\end{theorem}

\cref{thm:lower} follows from \cref{thm:lower-sys}, since a result of Rado \cite{Rad33} shows that the  linear equation $a_1x_1+\dots+a_kx_k=0$ has a monochromatic solution in every $2$-coloring of $[n]$, for $n$ sufficiently large. \cref{thm:lower-sys} strengthens an earlier result of Frankl, Graham, and R\"odl~\cite[Theorem 1]{FGR89}. They showed that if $Bx=0$ has a monochromatic solution in every coloring of $\NN$ with finitely many colors, then it has $\Omega_{B,r}(n^{k-w})$ monochromatic solutions in any $r$-coloring of $[n]$. In other words, they get the same conclusion as \cref{thm:lower-sys} under a stronger hypothesis.
\subsection*{Paper organization}
We set up some notation and definitions in \cref{sec:prelim}. We prove \cref{thm:lower,thm:lower-sys} in \cref{sec:lower} and \cref{thm:even-uncommon} in \cref{sec:even-uncommon}. The first results follow from the arithmetic removal lemma combined with a short supersaturation argument. The second result is proved by first showing that the commonness of this equation is equivalent to the eventual positive semidefiniteness of a sequence of matrices. We then provide an explicit counterexample to this positive semidefiniteness.

\subsection*{Acknowledgements}
The authors thank Shengtong Zhang for corrections on this manuscript.

\section{Preliminaries}\label{sec:prelim}

\subsection*{Notation}
Throughout, for $n \in \NN$, we let $[n] = \{1, 2, \ldots, n\}$.
We sometimes write $(x)_+$ for $\max(x,0)$ for $x \in \RR$.

let $L(x_1,\ldots,x_k) = a_1 x_1 + \cdots +a_k x_k$ denote a $k$-variable linear form. We will always assume that the coefficients are non-zero and at least one is positive and one is negative. For a function $f\colon [n]\to[0,1]$, define
  \[   T_L(f)=\sum_{\substack{x_1,\dots,x_k\in[n]:\\a_1x_1+\dots+a_kx_k=0}}f(x_1)\cdots f(x_k).
    \]
    This is the weighted count of solutions to $L=0$ in $f$. Also define
    \[
    t_L(f)= \frac{T_L(f)}{T_L(1)},
    \]
    the weighted density of solutions to $L=0$ in $f$.
    
    Throughout, we will assume that $n$ is sufficiently large. Together with the assumption on the signs of the coefficients of $L$, this implies that $T_L(1)>0$.
    
\subsection*{Common and locally common linear equations}
We first give precise definitions of {commonness} and {local commonness} of a linear equation. We make this definition for 2-colorings, though one could study the analogous notion of $r$-commonness for any $r$.
\begin{definition}\label{def:common}
    Let $L(x_1,\ldots,x_k)= a_1x_1+\dots+a_kx_k$ be a linear form. 
    \begin{enumerate}
        \item We say that $L =0$ is \textit{common over the integers} if for all sufficiently large $n$ and all functions $f\colon [n] \to [0, 1]$, we have $t_L(f)+t_L(1-f) \ge 2^{1-k}$.
        \item We say that $L=0$ is \textit{locally common over the integers} if for all sufficiently large $n$, there exists $\epsilon>0$ such that for all functions $f\colon [n]\to[-1,1]$ we have $t_L(1/2+\epsilon f)+t_L(1/2-\epsilon f) \ge 2^{1-k}$.
        \item We say that $L=0$ is \textit{weakly locally common over the integers} if for all sufficiently large $n$ and all functions $f\colon [n]\to\RR$ there exists $\epsilon_0(f)>0$ such that for all $\epsilon\in(0,\epsilon_0)$, we have $t_L(1/2+\epsilon f)+t_L(1/2-\epsilon f) \ge 2^{1-k}$.
    \end{enumerate}
\end{definition}
Note that if $A$ is a uniformly random 2-coloring of $[n]$, then $\EE[t_L(\mathbf 1_A)+t_L(1-\mathbf 1_A)]=2^{1-k}+o_{n\to\infty}(1)$. Thus $L=0$ is common if and only if the uniformly random 2-coloring asymptotically minimizes the number of monochromatic solutions.

\subsection*{Fourier transform on cyclic groups} The proof of \cref{thm:even-uncommon} utilizes basic Fourier analysis over $\ZZ/N\ZZ$, which we outline here. For every $N\in\NN$ and the cyclic group $\ZZ/N\ZZ$, let $\widehat{\ZZ/N\ZZ}$ denote the group of characters of $\ZZ/N\ZZ$, i.e., homomorphisms $\chi\colon\ZZ/N\ZZ\to\CC$. The Fourier transform of a function $f\colon\ZZ/N\ZZ\to\CC$ is defined by $$\widehat f(\chi)=\EE_{x\in \ZZ/N\ZZ}f(x)\overline{\chi(x)},$$ with the inverse Fourier transform $$f(x)=\sum_{\chi\in \widehat{\ZZ/N\ZZ}}\widehat f(\chi)\chi(x).$$
We will also use the simple fact that for all $x\in \ZZ/N\ZZ$, we have
\begin{align*}
    \sum_{\chi\in \widehat{\ZZ/N\ZZ}}\chi(x)=\begin{cases}
        N \quad & \text{if }x=0,\\
        0 & \text{if }x\neq 0.
    \end{cases}
\end{align*}

\section{Finding many monochromatic solutions}
\label{sec:lower}

In this section we show \cref{thm:lower-sys} (which implies \cref{thm:lower} by the discussion after the statement of \cref{thm:lower-sys}). This result follows from a type of supersaturation result when combined with the arithmetic removal lemma.

\begin{proposition}
\label{prop:removal-1}
For $B\in\ZZ^{w\times k}$ and $r\geq 2$, suppose that for $n$ sufficiently large, in every $r$-coloring of $[n]$ there is a monochromatic solution to the system $Bx=0$. Then there exists $\epsilon>0$ such that in every $r$-coloring of $[n]$ with $n$ sufficiently large, there exists a monochromatic solution to $Bx=0$ in every subset $A\subseteq[n]$ of size $|A|\geq (1-\epsilon)n$.
\end{proposition}

\begin{proof}
By hypothesis, there exists $N_0\geq 1$ so that every 2-coloring of $[N_0]$ contains a monochromatic solution to $Bx=0$. Observe that $[n]$ contains the subsets $S_r:=\{r,2r,\dots,N_0r\}$ for $1\leq r\leq \floor{n/N_0}$ and each $S_r$ contains some monochromatic solution $\mathbf x_r:=(x_{r,1},\dots,x_{r,k})$ to the equation $Bx=0$. 

Note that every $i\in[n]$ lies in at most $N_0$ different $S_r$'s. 
Thus, any $A\subseteq [n]$ of size $|A|>n-N_0^{-1}\floor{n/N_0}$ contains an entire set $S_r$ for some $r$ and thus contains the monochromatic solution $\mathbf x_r$. This implies the desired result with $\epsilon=N_0^{-2}/2$.
\end{proof}

\begin{proof}[Proof of \cref{thm:lower-sys}]
Consider a matrix $B\in\ZZ^{w\times k}$ and $r\geq 2$ as in the theorem statement. Take $\epsilon>0$ given by \cref{prop:removal-1}. By the arithmetic removal lemma \cite[Theorem 2]{KSV}, there exists $\delta=\delta(B,\epsilon/r)>0$ such that the following holds. For any $n$ and $X\subseteq[n]$ such that $Bx=0$ has fewer than $\delta n^{k-w}$ solutions in $X^k$, there exists a set $R\subseteq [n]$ of size at most $(\epsilon/r)n$ such that $X\setminus R$ has no solutions to $Bx=0$. 

For contradiction, suppose for sufficiently large $n$ there exists an $r$-coloring $[n]=X_1\sqcup \cdots\sqcup X_r$ that contains fewer than $\delta n^{k-w}$ monochromatic solutions to $Bx=0$. Then by the arithmetic removal lemma, we can remove at most $\epsilon n$ elements from $[n]$ to destroy all the monochromatic solutions. This violates \cref{prop:removal-1}.
\end{proof}

\section{Uncommon symmetric 4-variable equations}
\label{sec:even-uncommon}

In this section we study commonness of $4$-term linear equations of the form $ax_1+b x_2-ax_3-bx_4=0$. We show that this family of equations has some particularly nice properties. In particular, we show that this equation is common if and only if it is locally common. (Clearly any common linear equation is locally common, but the converse does not need to hold.) We also show that the equation $ax_1+b x_2-ax_3-bx_4=0$ is locally common (and thus common) if and only if a certain sequence of $n\times n$ matrices is eventually positive semidefinite.

\begin{proposition}
\label{thm:quadratic-quartic-decomp}
For $1\leq a<b$ coprime, define $L(x_1,x_2,x_3,x_4)=ax_1+b x_2-ax_3-b x_4$. There exist a quadratic form $\xi_{a,b}$ and a quartic form $\zeta_{a,b}$ such that
\[t_L(1/2+f)+t_L(1/2-f)=2^{1-k}+\xi_{a,b}(f)+\zeta_{a,b}(f)\]
for all $f\colon[n]\to\RR$. Furthermore, for all $f\colon[n]\to\RR$
\[\zeta_{a,b}(f)\geq 0.\]
\end{proposition}

\begin{proof}
Let $N=(a+b)n$. Then if $x_1,x_2,x_3,x_4\in[n]$ satisfy $ax_1+bx_2-ax_3-bx_4\equiv0\pmod N$, we see that $ax_1+bx_2-ax_3-bx_4=0$. For any function $f\colon[n]\to\RR$, extend $f$ to a function $f\colon\ZZ/N\ZZ\to\RR$ taking value 0 outside $[n]$. Then a standard calculation in Fourier analysis gives
\begin{align*}
T_L(f)&=\sum_{\substack{x_1,x_2,x_3,x_4\in\ZZ/N\ZZ:\\ax_1+b x_2-ax_3-b x_4=0}}f(x_1) f(x_2) f(x_3) f(x_4)\\
&=\frac1N\sum_{x_1,x_2,x_3,x_4\in\ZZ/N\ZZ}f(x_1) f(x_2) f(x_3) f(x_4)\sum_{\chi\in\widehat {\ZZ/N\ZZ}}\chi(ax_1+b x_2-ax_3-b x_4)\\
&=N^3\sum_{\chi\in\widehat {\ZZ/N\ZZ}}\EE_{x_1}[f(x_1)\chi(ax_1)]\EE_{x_2}[f(x_2)\chi(bx_2)]\EE_{x_3}[f(x_3)\chi(-ax_3)]\EE_{x_4}[f(x_4)\chi(-bx_4)]\\
&=N^3\sum_{\chi\in\widehat {\ZZ/N\ZZ}}|\widehat f(a\chi)|^2|\widehat f(b\chi)|^2.
\end{align*}
Let $I\colon\ZZ/N\ZZ\to \{0,1\}$ be the indicator function of $[n]$. Then,
\begin{align*}
T_L(1/2+f)+&T_L(1/2-f)-2T_L(1/2)\\
&= N^3\sum_{\chi\in\widehat {\ZZ/N\ZZ}}\left|\frac12 \widehat I(a\chi)+\widehat f(a\chi)\right|^2\left|\frac12 \widehat I(b\chi)+\widehat f(b\chi)\right|^2\\
&\qquad \qquad \qquad +\left|\frac12 \widehat I(a\chi)-\widehat f(a\chi)\right|^2\left|\frac12 \widehat I(b\chi)-\widehat f(b\chi)\right|^2-2\left|\frac12\widehat I(a\chi)\right|^2\left|\frac12\widehat I(b\chi)\right|^2\\
&= \Xi_{a,b}(f)+Z_{a,b}(f),
\end{align*}
where
\begin{align*}
\Xi_{a,b}(f)
&=\frac{N^3}2\sum_{\chi\in\widehat {\ZZ/N\ZZ}}\widehat f(b\chi){\widehat f(-b\chi)}\widehat I(a\chi){\widehat I(-a\chi)}
+{\widehat f(-a\chi)\widehat f(-b\chi)}\widehat I(a\chi)\widehat I(b\chi)\\
&\qquad+\widehat f(a\chi){\widehat f(-b\chi)\widehat I(-a\chi)}{\widehat I(b\chi)}
+{\widehat f(-a\chi)}{\widehat f(b\chi)\widehat I(a\chi)}{\widehat I(-b\chi)}\\
&\qquad+{\widehat f(a\chi)}{\widehat f(b\chi)}{\widehat I(-a\chi)\widehat I(-b\chi)}
+\widehat f(a\chi){\widehat f(-a\chi)}\widehat I(b\chi){\widehat I(-b\chi)},\\
Z_{a,b}(f)&=2N^3\sum_{\chi\in\widehat {\ZZ/N\ZZ}}\widehat f(a\chi)\widehat f(-a\chi)\widehat f(b\chi)\widehat f(-b\chi)=2N^3\sum_{\chi\in\widehat {\ZZ/N\ZZ}}|\widehat f(a\chi)|^2|\widehat f(b\chi)|^2.
\end{align*}
From this formula, we see that $$t_L(1/2+f)+t_L(1/2-f)=2t_L(1/2)+\xi_{a,b}(f)+\zeta_{a,b}(f)=2^{1-k}+\xi_{a,b}(f)+\zeta_{a,b}(f),$$ where $\xi_{a,b}(f)=\Xi_{a,b}(f)/T_L(1)$ and $\zeta_{a,b}(f)=Z_{a,b}(f)/T_L(1)$. In particular, $\xi_{a,b}$ is a quadratic form and $\zeta_{a,b}$ is a quartic form. Furthermore, since $Z_{a,b}(f)$ can be written as a sum of squares, we see that $\zeta_{a,b}(f)\geq 0$ for all $f$.
\end{proof}

Since $\xi_{a,b}(f)$ is a quadratic form, there exists a unique symmetric $n\times n$ matrix $m^{(a,b)}=m^{(a,b)}_n$ such that $\xi_{a,b}(f)=f^\intercal m_n^{(a,b)} f$ for all $f\colon[n]\to\RR$ (where we alternatively view such a function as a column vector $f\in\RR^n$).

\begin{corollary}
\label{thm:common-equiv}
For $1\leq a<b$ coprime, consider $L(x_1,x_2,x_3,x_4)=ax_1+bx_2-ax_3-bx_4$. The following are equivalent:
\begin{enumerate}
    \item $L=0$ is common over the integers;
    \item $L=0$ is locally common over the integers;
    \item $L=0$ is weakly locally common over the integers;
    \item The matrices $m^{(a,b)}_n$ are positive semidefinite for all sufficiently large $n$.
\end{enumerate}
\end{corollary}

\begin{proof}
Clearly (1) implies (2) implies (3). Furthermore, (3) implies (4), since (3) implies that for all sufficiently large $n$ and all $f\colon[n]\to\RR$, we have $\epsilon^2 f^\intercal m^{(a,b)}_n f+\epsilon^4 \zeta_{a,b}(f)\geq 0$ for all sufficiently small $\epsilon>0$. This implies that $f^\intercal m^{(a,b)}_n f\geq 0$ for all sufficiently large $n$, as desired. Finally, to show that (4) implies (1), let $f\colon[n]\to[-1/2,1/2]$ be any function. Then
\[t_L(1/2+f)+t_L(1/2-f)=2^{1-k}+f^\intercal m^{(a,b)}_n f+\zeta_{a,b}(f)\geq 2^{1-k},\]
by assumption on the positive semidefiniteness of $m^{(a,b)}_n$ and since \cref{thm:quadratic-quartic-decomp} implies that $\zeta_{a,b}$ is non-negative.
\end{proof}

To study the sequence $m^{(a,b)}_n$ of $n\times n$ matrices, we define a symmetric kernel $H_{a,b}\colon[0,1]^2\to\RR$ and then show that, informally, $H_{a,b}$ is a rescaled and ``smoothed'' version of the limit of this sequence of matrices. (We will see that the $(x,y)$ entry of $m^{(a,b)}_n$ depends on the value of $x-y\pmod{ab}$; in the definition of $H_{a,b}$, we will average over these contributions.)

Define the symmetric kernel $H_{a,b}\colon[0,1]^2\to\RR$ by
\begin{align*}
H_{a,b}(u,v)&=\alpha_{a,b}(au+bv)+\alpha_{a,b}(av+bu)+\alpha_{a,b}(a(1-v)+bu)+\alpha_{a,b}(a(1-u)+bv)\\&+b^{-1}\beta(a(u-v)/b )+a^{-1}\beta(b (u-v)/a),
\end{align*}
where the auxiliary functions $\alpha_{a,b},\beta \colon\RR\to\RR_{\geq 0}$ are defined by
\begin{equation}\label{e:alpha}
\alpha_{a,b}(u)=\begin{cases}
\frac{u}{ab} \quad&\text{if }0\leq u\leq a,\\
\frac 1b  &\text{if }a\leq u\leq b,\\
\frac{a+b-u}{ab} &\text{if }b\leq u\leq a+b,\\
0&\text{otherwise},
\end{cases}
\end{equation}
and
\begin{equation}\label{e:beta}
    \beta(u)=\begin{cases}
u+1\quad &\text{if }-1\leq u\leq 0,\\
1-u&\text{if }0\leq u\leq 1,\\
0&\text{otherwise}.
\end{cases}
\end{equation}

\begin{proposition}
\label{thm:non-psd-kernel}
For $1\leq a<b$ relatively prime, if $H_{a,b}$ is not positive semidefinite, then $m^{(a,b)}_n$ is not positive semidefinite for infinitely many $n$.
\end{proposition}

\begin{proof}
In the proof of \cref{thm:quadratic-quartic-decomp}, we computed $m^{(a,b)}_n$ in frequency space. To prove this result, we will compute these matrices in physical space. For $f\colon [n]\to\RR$ we can write
\begin{align*}
T_L(1/2+ f)+T_L(1/2- f)&=
\sum_{\substack{x,x',y,y'\in[n]:\\ax + b y = ax' + b y'}}\left( \frac12 +  f(x) \right)  \left( \frac12 +  f(y)\right)  \left( \frac12 +  f(x') \right)  \left( \frac12 +  f(y')\right) \\
&\qquad+ \left( \frac12 -  f(x)\right)  \left( \frac12 -  f(y) \right)  \left( \frac12 -  f(x') \right)  \left( \frac12 -  f(y') \right) \\
&=\sum_{\substack{x,x',y,y'\in[n]:\\ax + b y = ax' + b y'}}\frac18+\frac12\big( f(x)f(y)+f(x)f(x')+f(x)f(y')\\
&\qquad+f(y)f(x')+f(y)f(y')+f(x')f(y')\big)+2f(x)f(y)f(x')f(y')\\
&= \frac18T_L(1) + T_L(1) \xi_{a,b}(f) +T_L(1) \zeta_{a,b}(f),
\end{align*}
where $\xi_{a,b}(f)$ and $\zeta_{a,b}(f)$ are the quadratic form and quartic form given in \cref{thm:quadratic-quartic-decomp}. The quadratic form corresponds to the middle six terms in the summation. Expanding, we write
\begin{align*}
2T_L(1)\xi_{a,b}(f) &=\sum_{\substack{x,x',y,y'\in[n]: \\ax + b y = ax' + b y'}} f(x)f(y)+f(x)f(x')+f(x)f(y')\\
&\qquad+f(y)f(x')+f(y)f(y')+f(x')f(y')\\
&=\sum_{x,y\in[n]} f(x)f(y)|\{(x',y')\in[n]^2:ax+by=ax'+by'\}|\\
&\qquad+\sum_{x,x'\in[n]}f(x)f(x')|\{(y,y')\in[n]^2:ax+by=ax'+by'\}|\\
&\qquad+\sum_{x,y'\in[n]}f(x)f(y')|\{(x',y)\in[n]^2:ax+by=ax'+by'\}|\\
&\qquad+\sum_{y,x'\in[n]}f(y)f(x')|\{(x,y')\in[n]^2:ax+by=ax'+by'\}|\\
&\qquad+\sum_{y,y'\in[n]}f(y)f(y')|\{(x,x')\in[n]^2:ax+by=ax'+by'\}|\\
&\qquad+\sum_{x',y'\in[n]}f(x')f(y')|\{(x,y)\in[n]^2:ax+by=ax'+by'\}|.
\end{align*}
Each of the expressions is easy to explicitly compute up to an error of $O(1)$. Recall that $(x)_+ = \max(x,0)$. The computation gives
\begin{align*}
2&T_L(1)\xi_{a,b}(f) =\sum_{x,y\in[n]} f(x)f(y)\paren{\frac1a\paren{\min\paren{n,\frac{ax+by-a}b}-\max\paren{\frac{ax+by-an}b,1}}_{+}+O(1)}\\
&\qquad+\sum_{x,x'\in[n]}f(x)f(x')\left(\left(\min\paren{n,\frac{ax-ax'+bn}b}\right.\right.\\&\qquad\qquad\qquad\qquad\left.\left.-\max\paren{\frac{ax-ax'+b}b,1}\right)_{+}+O(1)\right)\cdot \mathbf 1(x\equiv x'\pmod b)\\
&\qquad+\sum_{x,y'\in[n]}f(x)f(y')\paren{\frac1a\paren{\min\paren{n,\frac{by'-ax+an}b}-\max\paren{\frac{by'-ax+a}b,1}}_{+}+O(1)}\\
&\qquad+\sum_{y,x'\in[n]}f(y)f(x')\paren{\frac1a\paren{\min\paren{n,\frac{by-ax'+an}b}-\max\paren{\frac{by-ax'+a}b,1}}_{+}+O(1)}\\
&\qquad+\sum_{y,y'\in[n]}f(y)f(y')\left(\left(\min\paren{n,\frac{by-by'+an}a}\right.\right.\\&\qquad\qquad\qquad\qquad\left.\left.-\max\paren{\frac{by-by'+a}a,1}\right)_{+}+O(1)\right)\cdot \mathbf 1(y\equiv y'\pmod a)\\
&\qquad+\sum_{x',y'\in[n]}f(x')f(y')\paren{\frac1a\paren{\min\paren{n,\frac{ax'+by'-a}b}-\max\paren{\frac{ax'+by'-an}b,1}}_{+}+O(1)}\\
&=\sum_{x,y\in[n]}f(x)f(y)\left(\frac1a\paren{\min\paren{n,\frac{ax+by-a}b}-\max\paren{\frac{ax+by-an}b,1}}_{+}\right.\\
&\qquad+\paren{\min\paren{n,\frac{ax-ay+bn}b}-\max\paren{\frac{ax-ay+b}b,1}}_{+}\cdot \mathbf 1(x\equiv y\pmod b)\\
&\qquad+\frac1a\paren{\min\paren{n,\frac{bx-ay+an}b}-\max\paren{\frac{bx-ay+a}b,1}}_{+}\\
&\qquad+\frac1a\paren{\min\paren{n,\frac{bx-ay+an}b}-\max\paren{\frac{bx-ay+a}b,1}}_{+}\\
&\qquad+\paren{\min\paren{n,\frac{bx-by+an}a}-\max\paren{\frac{bx-by+a}a,1}}_{+}\cdot \mathbf 1(x\equiv y\pmod a)\\
&\qquad+\left.\frac1a\paren{\min\paren{n,\frac{ax+by-a}b}-\max\paren{\frac{ax+by-an}b,1}}_{+}+O(1)\right) \\
&= \sum_{x, y \in [n]} f(x) f(y) A(x, y),
\end{align*}
where we define $A(x,y)$ to be the expression in parentheses above.

Observe that the first and sixth terms are identical and the third and fourth terms are identical. Recalling that $\xi_{a,b}(f)=f^\intercal m^{(a,b)} f$ where $m$ is a symmetric $n\times n$ matrix, the above expansion implies that the entries of $m^{(a,b)}$ satisfy the equation $4T_L(1)m^{(a,b)}_{x,y}=A(x,y)+A(y,x)$. Let us expand this and also round the expression by $O(1)$ to simplify it slightly.
\begin{align*}
4T_L(1) m^{(a,b)}_{x,y}&= \frac2a\paren{\min\paren{n,\frac{ax+by}b}-\paren{\frac{ax+by-an}b}_+}_{+}\\
&\qquad+\frac2a\paren{\min\paren{n,\frac{ay+bx}b}-\paren{\frac{ay+bx-an}b}_+}_{+}\\
&\qquad+\paren{\min\paren{n,\frac{ax-ay+bn}b}-\paren{\frac{ax-ay}b}_+}_{+}\cdot \mathbf 1(x\equiv y\pmod b)\\
&\qquad+\paren{\min\paren{n,\frac{ay-ax+bn}b}-\paren{\frac{ay-ax}b}_+}_{+}\cdot \mathbf 1(x\equiv y\pmod b)\\
&\qquad+\frac2a\paren{\min\paren{n,\frac{bx-ay+an}b}-\paren{\frac{bx-ay}b}_+}_{+}\\
&\qquad+\frac2a\paren{\min\paren{n,\frac{by-ax+an}b}-\paren{\frac{by-ax}b}_+}_{+}\\
&\qquad+\paren{\min\paren{n,\frac{bx-by+an}a}-\paren{\frac{bx-by}a}_+}_{+}\cdot \mathbf 1(x\equiv y\pmod a)\\
&\qquad+\paren{\min\paren{n,\frac{by-bx+an}a}-\paren{\frac{by-bx}a}_+}_{+}\cdot \mathbf 1(x\equiv y\pmod a)+O(1)\\
&=\frac{2n}{ab}\paren{\min\paren{b,\frac{ax+by}n}-\paren{\frac{ax+by}n-a}_+}_{+}\\
&\qquad+\frac{2n}{ab}\paren{\min\paren{b,\frac{ay+bx}n}-\paren{\frac{ay+bx}n-a}_+}_{+}\\
&\qquad+n\paren{\min\paren{1,\frac{ax-ay}{bn}+1}-\paren{\frac{ax-ay}{bn}}_+}_{+}\cdot \mathbf 1(x\equiv y\pmod b)\\
&\qquad+n\paren{\min\paren{1,\frac{ay-ax}{bn}+1}-\paren{\frac{ay-ax}{bn}}_+}_{+}\cdot \mathbf 1(x\equiv y\pmod b)\\
&\qquad+\frac{2n}{ab}\paren{\min\paren{b,\frac{bx-ay}n+a}-\paren{\frac{bx-ay}n}_+}_{+}\\
&\qquad+\frac{2n}{ab}\paren{\min\paren{b,\frac{by-ax}n+a}-\paren{\frac{by-ax}n}_+}_{+}\\
&\qquad+n\paren{\min\paren{1,\frac{bx-by}{an}+1}-\paren{\frac{bx-by}{an}}_+}_{+}\cdot \mathbf 1(x\equiv y\pmod a)\\
&\qquad+n\paren{\min\paren{1,\frac{by-bx}{an}+1}-\paren{\frac{by-bx}{an}}_+}_{+}\cdot \mathbf 1(x\equiv y\pmod a)+O(1)\\
&:=B(x,y).
\end{align*}

Note that the auxiliary functions that we defined earlier in~\cref{e:alpha,e:beta} satisfy
\[\alpha_{a,b}(u)=\frac1{ab}\paren{\min(b,u)-(u-a)_+}_+\]
and
\[\beta(u)=\beta(-u)=(\min(1,u+1)-u_+)_+.\]
Now we can simplify the formula for $B$. For $u,v\in[0,1]$ such that $un,vn\in[n]$, we see that
\begin{align*}
\frac1{2n}&B(un,vn)=\alpha_{a,b}(au+bv)+\alpha_{a,b}(av+bu)+\beta(a(x-y)/b)\mathbf1(un\equiv vn\pmod b)\\
&+\alpha_{a,b}(a(1-v)+bu)+\alpha_{a,b}(a(1-u)+bv)+\beta(b (x-y)/a)\mathbf1(un\equiv vn\pmod a)+O(1/n).
\end{align*}
One should think of the indicator $\mathbf1(un\equiv vn\pmod b)$ having average value $1/b$ and $\mathbf1(un\equiv vn\pmod a)$ having average value $1/a$. Thus one can see that $H_{a,b}(u,v)$ is a ``smoothed'' approximation of $B(un,vn)/2n$. We use this intuition to prove the desired result.

Now if $H_{a,b}$ is not positive semidefinite, there exists some $\phi \in L^2([0,1])$ such that
\[\iint_{[0,1]^2}H_{a,b}(u,v)\phi(u)\phi(v)\,\text du\text dv=-\delta<0.\]
Let $n$ be a sufficiently large multiple of $ab$. Define $f_n\colon [n]\to\RR$ by
\[f_n(abi+j)=\frac1{ab}\int_{abi/n}^{ab(i+1)/n}\phi(u)\,\text du\]
for $j\in[ab]$ and $0\leq i<n/ab$.

We claim that $f_n^\intercal m^{(a,b)}_n f_n<0$. Since $\tfrac{2T_L(1)}{n}m^{(a,b)}_n=\tfrac1{2n}B$, it suffices to show that $f_n^\intercal \tfrac1{2n}B f_n<0$. Now note that for fixed $0\leq i,i'<n/ab$,
\begin{align*}
\sum_{j,j'\in[ab]}& f_n(abi+j)f_n(abi'+j') \frac1{2n} B(abi+j,abi'+j')\\
&= f_n(ab(i+1))f_n(ab(i'+1))\sum_{j,j'\in[ab]} \frac1{2n} B(abi+j,abi'+j')\\
&= a^2b^2 f_n(ab(i+1))f_n(ab(i'+1))\paren{H_{a,b}(ab(i+1)/n,ab(i'+1)/n)+O(1/n)}\\
&= \iint_{(abi/n,abi'/n)+[0,ab/n]^2}(H_{a,b}(u,v)+O(1/n))\phi(u)\phi(v)\,\text du\text dv\\
&= \iint_{(abi/n,abi'/n)+[0,ab/n]^2}H_{a,b}(u,v)\phi(u)\phi(v)\,\text du\text dv\\
&\qquad+O\paren{\frac{1}{n} \int_{abi/n}^{ab(i+1)/n} |\phi(u)|\,\text du\int_{abi'/n}^{ab(i'+1)/n} |\phi(v)|\,\text dv}.
\end{align*}
The third line follows since $\alpha_{a,b},\beta$ are 1-Lipshitz. The fourth line follows since $H$ is $O(1)$-Lipshitz. Then summing over $i,i'$, we conclude that
\[f_n^\intercal \frac{B}{2n} f_n=\iint_{[0,1]^2}H_{a,b}(u,v)\phi(u)\phi(v)\,\text du\text dv+O\paren{\frac{\|\phi\|_{1}^2}{n}}=-\delta+O\paren{\frac{\|\phi\|_{1}^2}{n}}.\]
For every $n$ that is a multiple of $ab$ and that is sufficiently large in terms of $\|\phi\|_1\leq 1+\|\phi\|_2^2<\infty$, this expression is negative. This proves the desired result.
\end{proof}

We are now ready to prove \cref{thm:even-uncommon}. We do so by exhibiting a specific test function that demonstrates that $H_{1,2}$ is not positive semidefinite. This function was found with computer assistance. We discretized $H_{1,2}$ into a large matrix, computed the eigenvector corresponding to its least eigenvalue, and then rounded that to produce the test function $\phi$ defined below.

\begin{proof}[Proof of \cref{thm:even-uncommon}]
We show that $x_1+2x_2-x_3-2x_4=0$ is uncommon over the integers. By \cref{thm:common-equiv,thm:non-psd-kernel}, it suffices to show that $H_{1,2}$ is not positive semidefinite. Consider the following explicit test function $\phi\colon[0,1]\to\RR$:
\[
\phi(x)=
\begin{cases}
-2\qquad&\text{if }\frac{0}{200}\leq x<\frac{7}{200},\\
10&\text{if }\frac{7}{200}\leq x<\frac{16}{200},\\
-18&\text{if }\frac{16}{200}\leq x<\frac{33}{200},\\
15&\text{if }\frac{33}{200}\leq x<\frac{67}{200},\\
-19&\text{if }\frac{67}{200}\leq x<\frac{84}{200},\\
9&\text{if }\frac{84}{200}\leq x<\frac{93}{200},\\
-4&\text{if }\frac{93}{200}\leq x<\frac{107}{200},\\
9&\text{if }\frac{107}{200}\leq x<\frac{116}{200},\\
-19&\text{if }\frac{116}{200}\leq x<\frac{133}{200},\\
15&\text{if }\frac{133}{200}\leq x<\frac{167}{200},\\
-18&\text{if }\frac{167}{200}\leq x<\frac{184}{200},\\
10&\text{if }\frac{184}{200}\leq x<\frac{193}{200},\\
-2&\text{if }\frac{193}{200}\leq x\leq \frac{200}{200}.
\end{cases}
\]
By explicit computation,
\[\iint_{[0,1]^2}H_{1,2}(x,y)\phi(x)\phi(y)\,\text dx\text dy=-\frac{120959}{1600000}<0.\qedhere\]
\end{proof}

\section{Further remarks}

Let us now summarize the state of affairs on the commonness of a single equation. Consider the linear form $L(x_1,\ldots,x_k)=a_1x_1+\cdots+a_kx_k$ with non-zero coefficients. If $k$ is odd, it is plausible that $L=0$ is uncommon for every such $L$. 

\begin{conjecture}
For $k\geq 3$ odd and $a_1,\ldots,a_k\in\ZZ\setminus\{0\}$, the equation $a_1x_1+\cdots+a_kx_k=0$ is uncommon over the integers.
\end{conjecture}

The $k=3$ case was proved by Costello and Elvin. For most equations of odd length, their uncommonness follows from the work of Fox--Pham--Zhao \cite{FPZ21} over $\FF_p$ via the reduction given in the introduction. However, there are still many equations of odd length whose uncommonness is hard to resolve via this technique. One particularly difficult case seems to be the equation $x_1+\cdots+x_{\ell+1}-x_{\ell_2}-\cdots-x_{2\ell+1}=0$.

If $k$ is even, then $L=0$ is uncommon over the integers unless the coefficients can be partitioned into pairs summing to zero. This also follows from the reduction from $\FF_p$ and the work of Fox--Pham--Zhao. However, if $k$ is even and the coefficients can be partitioned into pairs summing to zero, we do not know if $L=0$ is common over the integers and we do not even have a plausible conjecture. We showed that the specific equation $x_1+2x_2-x_3-2x_4=0$ is not common over the integers. On the other hand, Costello and Elvin showed that the equation $x_1+\cdots+x_{k/2}-x_{k/2+1}-\cdots-x_{k}=0$ is common over the integers for each even $k$ \cite[Appendix A]{CE23}. To date, these are the only equations known to be common over the integers.

\begin{question}
For which $k$ even and $a_1,\ldots,a_{k/2}\in\ZZ\setminus\{0\}$, is the equation $a_1x_1+\cdots+a_{k/2}x_{k/2}-a_1x_{k/2+1}-\cdots-a_{k/2}x_k=0$ common over the integers?
\end{question}

Even the $k=4$ case remains elusive.

\begin{question}
For which pairs $(a,b)$, is the equation $ax_1+bx_2-ax_3-bx_4=0$ common over the integers?
\end{question}

In \cref{sec:even-uncommon}, we showed that every equation  $ax_1+b x_2-ax_3-b x_4=0$ is uncommon over the integers if the operator $H_{a,b}$ is not positive semidefinite. We are able to show that $x_1+2x_2-x_3-2x_4=0$ is uncommon by explicitly demonstrating that the operator $H_{1,2}$ is not positive semidefinite. Limited computational evidence suggests that for $1\le a<b$, the operator $H_{a,b}$ is not positive semidefinite when $b/a$ is small, but is positive semidefinite (but not positive definite) for $b/a$ large. For any single $(a,b)$ with $b/a$ small, one should be able to explicitly demonstrate a negative eigenvalue of $H_{a,b}$ as in the proof of \cref{thm:even-uncommon}. Proving this result uniformly for all $(a,b)$ with $b/a$ small seems to be more demanding. When $H_{a,b}$ is positive semidefinite, our argument does not yield a conclusion about commonness of the corresponding equation.

We believe that one particularly interesting regime is when $b$ is much larger than $a$, where computational evidence suggests that $H_{a,b}$ is positive semidefinite, but not positive definite.

\bibliographystyle{plain}
\bibliography{common.bib}

\begin{thebibliography}{10}

\bibitem{BCG10}
Steve Butler, Kevin~P. Costello, and Ron Graham.
\newblock Finding patterns avoiding many monochromatic constellations.
\newblock {\em Experiment. Math.}, 19(4):399--411, 2010.

\bibitem{CE23}
Kevin~P. Costello and Gabriel Elvin.
\newblock Avoiding monochromatic solutions to 3-term equations.
\newblock {\em J. Comb.}, 14(3):281--304, 2023.

\bibitem{Dat03}
Boris~A. Datskovsky.
\newblock On the number of monochromatic {S}chur triples.
\newblock {\em Adv. in Appl. Math.}, 31(1):193--198, 2003.

\bibitem{FPZ21}
Jacob Fox, Huy~Tuan Pham, and Yufei Zhao.
\newblock Common and {S}idorenko linear equations.
\newblock {\em Q. J. Math.}, 72(4):1223--1234, 2021.

\bibitem{FGR89}
P.~Frankl, R.~L. Graham, and V.~R\"odl.
\newblock On the distribution of monochromatic configurations.
\newblock In {\em Irregularities of partitions ({F}ert\H od, 1986)}, volume~8
  of {\em Algorithms Combin. Study Res. Texts}, pages 71--87. Springer, Berlin,
  1989.

\bibitem{GRR96}
Ronald Graham, Vojtech R\"odl, and Andrzej Ruci\'nski.
\newblock On {S}chur properties of random subsets of integers.
\newblock {\em J. Number Theory}, 61(2):388--408, 1996.

\bibitem{GKW24}
Torin Greenwood, Jonathan Kariv, and Noah Williams.
\newblock From discrete to continuous: monochromatic 3-term arithmetic
  progressions.
\newblock {\em Math. Comp.}, 93(350):2959--2983, 2024.

\bibitem{KSV}
Daniel Kr{\'a}l', Oriol Serra, and Llu{\'\i}s Vena.
\newblock A removal lemma for systems of linear equations over finite fields.
\newblock {\em Israel Journal of Mathematics}, 187:193--207, 2012.

\bibitem{LP12}
Linyuan Lu and Xing Peng.
\newblock Monochromatic 4-term arithmetic progressions in 2-colorings of
  {$\mathbb{Z}_n$}.
\newblock {\em J. Combin. Theory Ser. A}, 119(5):1048--1065, 2012.

\bibitem{PRS08}
Pablo~A. Parrilo, Aaron Robertson, and Dan Saracino.
\newblock On the asymptotic minimum number of monochromatic 3-term arithmetic
  progressions.
\newblock {\em J. Combin. Theory Ser. A}, 115(1):185--192, 2008.

\bibitem{Rad33}
Richard Rado.
\newblock Studien zur {K}ombinatorik.
\newblock {\em Math. Z.}, 36(1):424--470, 1933.

\bibitem{RZ98}
Aaron Robertson and Doron Zeilberger.
\newblock A {$2$}-coloring of {$[1,n]$} can have {$(1/22)n^2+O(n)$}
  monochromatic {S}chur triples, but not less!
\newblock {\em Electron. J. Combin.}, 5:Research Paper 19, 4, 1998.

\bibitem{Sch99}
Tomasz Schoen.
\newblock The number of monochromatic {S}chur triples.
\newblock {\em European J. Combin.}, 20(8):855--866, 1999.

\bibitem{TW17}
Thotsaporn Thanatipanonda and Elaine Wong.
\newblock On the minimum number of monochromatic generalized {S}chur triples.
\newblock {\em Electron. J. Combin.}, 24(2):Paper No. 2.20, 20, 2017.

\end{thebibliography}

\end{document}